\documentclass[12pt]{article}
\usepackage{amssymb}
\usepackage{amsthm}
\usepackage{amsmath}    
\usepackage{pgf,tikz}
\usepackage{graphicx}
\usepackage{amssymb}
\usepackage{graphics}
\usepackage{amsthm}

\input pictex.tex

\theoremstyle{definition}

\newtheorem{thm}{Theorem}[section]

\newtheorem{lem}[thm]{Lemma}

\newtheorem{rem}[thm]{Remark}

\newcommand{\Diff}{\mathrm{Diff}}
\newcommand{\var}{\mathrm{var}}

\begin{document}

\date{}
\author{Leonardo Dinamarca \quad  \& \quad Maximiliano Escayola}

\title{Some examples of distorted interval diffeomorphisms\\ 
of intermediate regularity}
\maketitle

\vspace{-0.4cm}

\noindent {\bf Abstract:} We improve a recent construction of Andr\'es Navas to produce the first examples of $C^2$-undistorted diffeomorphisms 
of the interval that are  $C^{1+\alpha}$-distorted (for every $\alpha < 1$).  We do this via explicit computations due to the failure of an extension 
to class $C^{1+\alpha}$ of a classical lemma related to the work of Nancy Kopell.

\vspace{0.2cm}

\noindent{\bf Keywords:} Interval diffeomorphism, H\"older continuity, distorted element, distortion.

\vspace{0.2cm}

\noindent{\bf 2010 Mathematics Subject Classification:} 20F99, 37C05, 37C85, 37E05. 

\vspace{0.7cm}

\noindent{\bf {\Large Introduction}}

\vspace{0.5cm}

We start by recalling the terminology introduced by Michail Gromov \cite{gromov}. Given a finitely generated group $\Gamma$, we fix a finite system 
of generators, and we denote $\| \cdot \|$ the corresponding word-length. An element $f \in \Gamma$ is said to be distorted if 
$$\lim_{n \to \infty} \frac{\| f^n \|}{n} = 0.$$
(Notice that this condition does not depend on the choice of the finite generating system.) Given an arbitrary group $G$, an element 
$f \in G$ is said to be distorted if there exists a finitely generated subgroup $\Gamma \subset G$ containing $f$ so that 
$f$ is distorted in $\Gamma$ in the sense above. 

Examples of ``large'' groups for which this notion becomes interesting are groups of diffeomorphisms of compact 
manifolds $M$. Very little is known about distorted elements therein. In particular, the following question from \cite{Na} is  
widely open: Given $r < s$, does there exist an undistorted element $f \in \Diff^{s}_+ (M)$ that is distorted when considered 
as an element of $\Diff^{r}_+ (M)$~? In \cite{Na}, Andr\'es Navas proves that this is the case for $M$ the closed interval, $r = 1$ and $s=2$. 
Actually, undistortion holds in the larger group $\Diff^{1+bv}_+([0,1])$ of $C^1$ diffeomorphisms with derivative of bounded variation. 

In this Note, we give an extension of this result from $C^1$ to $C^{1+\alpha}$ regularity.

\vspace{0.5cm}

\noindent{\bf Theorem.} {\em There exist $C^{\infty}$ diffeomorphisms of $[0,1]$ that are distorted in 
$\Diff^{1+\alpha}_+([0,1])$ for all $\alpha > 0$ yet undistorted in $\Diff^{1+bv}_+([0,1])$.}

\vspace{0.5cm}

The groups we consider are variations of those introduced in \cite{Na}. One of the new contribution consists in improving the regularity 
of some elements, which is not at all straightforward. Indeed, the construction of \cite{Na} uses a well-known lemma that ensures $C^1$ 
regularity of maps built by pasting together infinitely many diffeomorphisms that are defined on disjoint intervals and satisfy certain equivariance 
relations. This  idea comes from the thesis of Nancy Kopell \cite{Ko}, and has been systematically used in the study of codimension-1 
foliations \cite{DF} and centralizers of diffeomorphisms \cite{BF}. 
Nevertheless, such a lemma is unavailable in $C^{1+\alpha}$ regularity and, as we show in the Appendix, 
it cannot hold without imposing extra hypothesis. We are hence forced to go into more explicit constructions and very long 
computations, which are however interesting by themselves. To do this, we use a classical technique of Dennis Pixton (later 
extended by Takashi Tsuboi \cite{tsuboi}) to produce commuting diffeomorphisms and control their $C^{1+\alpha}$ norms.

Here and in what follows, all maps we consider are supposed to preserve the orientation.


\section{On a family of $C^{1+bv}$-undistorted diffeomorphisms}

Recall that for a $C^{1+bv}$ diffeomorphism $f$ of a compact 1-manifold, the {\em asymptotic distortion} was 
defined by Navas in \cite{Na-or} as 
$$\mathrm{dist}_{\infty} (f) := \lim_{n \to \infty} \frac{\mathrm{var} (\log Df^n)}{n}.$$
By the subaditivity of $\var (\log D(\cdot))$, if $f$ is a distorted element of the group of $C^{1+bv}$ diffeomorphisms, 
then $\mathrm{dist}_{\infty} (f) = 0$.

The family of diffeomorphisms with positive asymptotic distortion studied in \cite{Na} is as follows: Start with a $C^{1+bv}$ diffeomorphism of 
$[0,1]$ with vanishing asymptotic distortion and no fixed point in $]0,1[$. Let $I$ be a {\em fundamental domain} for the action of $f$, that is, 
an open interval with endpoints $x_0$ and $x_1 := f (x_0)$ for a certain $x_0 \in \, ]0,1[$. Let $g$ be any nontrivial $C^{1+bv}$ 
diffeomorphism of $]0,1[$ supported on $I$. Then the diffeomorphism $\bar{f} := fg$ has positive asymptotic distortion and, 
in particular, it is undistorted in $\mathrm{Diff}_+^{1+bv} ([0,1])$ (hence in $\mathrm{Diff}_+^2([0,1])$). 
This fact follows from \cite{EN} (see Lemmas 2.2 and 7.2 therein) by using the relation between the asymptotic distortion and the Mather 
invariant. For the reader's convenience, below we present a short and direct argument based on Kopell's like estimates \cite{Ko,book}.
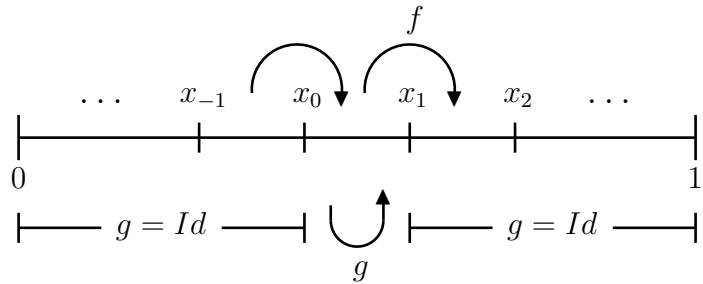
\vspace{-12pt}
\begin{center}
\begin{tikzpicture}[x=1.0cm,y=1.0cm]
\draw [line width=1.pt] (0.,2.3)-- (0.,1.7);
\draw [line width=1.pt] (0.,2.)-- (9.,2.);
\draw [line width=1.pt] (9.,2.3)-- (9.,1.7);
\draw [line width=1.pt] (2.4,2.2)-- (2.4,1.8);
\draw [line width=1.pt] (3.8,2.2)-- (3.8,1.8);
\draw [line width=1.pt] (5.2,2.2)-- (5.2,1.8);
\draw [line width=1.pt] (6.6,2.2)-- (6.6,1.8);
\draw [line width=1.pt] (0.,1.0)-- (0.,0.6);
\draw [line width=1.pt] (3.8,1.0)-- (3.8,0.6);
\draw [line width=1.pt] (5.2,1.0)-- (5.2,0.6);
\draw [line width=1.pt] (9.,1.0)-- (9.,0.6);
\draw [line width=1.pt] (5.2,0.8)-- (6.3,0.8);
\draw [line width=1.pt] (7.9,0.8)-- (9.,0.8);
\draw [line width=1.pt] (2.7,0.8)-- (3.8,0.8);
\draw [line width=1.pt] (0.,0.8)-- (1.1,0.8);
\draw (2.0,2.8) node[anchor=north west] {$x_{-1}$};
\draw (3.5,2.8) node[anchor=north west] {$x_{0}$};
\draw (4.9,2.8) node[anchor=north west] {$x_1$};
\draw (6.3,2.8) node[anchor=north west] {$x_{2}$};
\draw (7.41,2.65) node[anchor=north west] {\scalebox{1.2}{$.\,.\,.$}};
\draw (0.656,2.65) node[anchor=north west] {\scalebox{1.2}{$.\,.\,.$}};
\draw (1.15,1.15) node[anchor=north west] {$g=Id$};
\draw (6.35,1.15) node[anchor=north west] {$g=Id$};
\draw (4.3,0.5) node[anchor=north west] {$g$};
\draw (-0.25,1.75) node[anchor=north west] {$0$};
\draw (8.75,1.75) node[anchor=north west] {$1$};
\draw (5.0,3.9) node[anchor=north west] {$f$};
\draw [shift={(3.7,2.64)},line width=1.pt]  plot[domain=-0.07:3.22,variable=\t]({1.*0.60*cos(\t r)+0.*0.60*sin(\t r)},{0.*0.60*cos(\t r)+1.*0.60*sin(\t r)});
\draw [shift={(5.2,2.64)},line width=1.pt]  plot[domain=-0.07:3.22,variable=\t]({1.*0.60*cos(\t r)+0.*0.60*sin(\t r)},{0.*0.60*cos(\t r)+1.*0.60*sin(\t r)});
\draw [shift={(4.5,0.9)},line width=1.pt]  plot[domain=-3.14:0.,variable=\t]({1.*0.35*cos(\t r)+0.*0.35*sin(\t r)},{0.*0.35*cos(\t r)+1.*0.35*sin(\t r)});
\draw [line width=1.pt] (4.85,0.9)-- (4.85,1.2);
\draw [line width=1.pt] (4.15,0.9)-- (4.15,1.1);
\draw (4.58,1.45) node[anchor=north west] {\scalebox{.8}{$\blacktriangle$}};
\draw (4.025,2.75) node[anchor=north west] {\scalebox{.8}{$\blacktriangledown$}};
\draw (5.525,2.75) node[anchor=north west] {\scalebox{.8}{$\blacktriangledown$}};
\end{tikzpicture}
\end{center}

 Let us consider the product $\bar{f}^n f^{-n}$. Since $f$ has vanishing asymptotic distortion, if we show that 
 $\var (\log D (\bar{f}^n f^{-n}) )$ has linear growth, the same will hold for $\var (\log D \bar{f}^n)$. Now, notice that 
 $$\bar{f}^n f^{-n} = (fgf^{-1}) (f^2 g f^{-2}) (f^3 gf^{-3}) \cdots (f^{n} g f^{-n})$$
 has support in the union of the intervals $f(I), f^2(I), \ldots, f^n (I)$, and equals $f^k g f^{-k}$ on each such interval $f^k (I)$.  
 In particular, its derivative at the endpoints $x_k,x_{k+1}$ of each of these intervals equals 1. We claim that there is a constant $\lambda > 1$ 
 such that, for all $k \geq 1$, there is a point $y_k \in f^k (I)$ satisfying $D (f^k g f^{-k}) (y_k) \geq \lambda$. Assuming this, we conclude 
 $$\var (\log D (\bar{f}^n f^{-n}) ) \geq \sum_{k=1}^{n} \big| \log D(f^kgf^{-k}) (y_k) - \log D (f^kgf^{-k}) (x_k) \big| \geq n \, \log (\lambda),$$
 which yields the desired linear growth.
 
 Now, to check the existence of $\lambda$ and the points $y_k$, let $V := \var (\log Df)$, and let $N \geq 1$ be such that $Dg^N (z) > e^{2V}$ 
 holds for some $z \!\in\! I$. We claim that $\lambda := e^{V/N}$ works. Assume otherwise. Then, for a certain $K \geq 0$, one would 
 have $ \| D (f^K g f^{-K} ) \|_{\infty} \leq e^{V / N}$, which by the chain rule would yield $ \| D (f^K g^N f^{-K} ) \|_{\infty} \leq e^{V}$. 
 However, at the point $z_K \!:=\! f^K (z)$, we have $D (f^K g^N f^{-K}) (z_K) > e^V$. Indeed,
\begin{eqnarray*}
\log (D (f^K g^N f^{-K}) (z_K)) 
&=& \log Df^K (g^N (z) ) + \log Dg^N (z) - \log Df^K (z) \\
&\geq& \log Dg^N (z)  - \big| \log Df^K (g^N (z) ) - \log Df^K (z) \big| \\
&>& 2V - \sum^{K-1}_{k=0} \big| \log Df (f^k (g^N (z))) - \log Df (f^k (z)) \big|. 
\end{eqnarray*}
Since both $z$ and $g^N (z)$ lie in the fundamental domain $I$ of $f$, 
$$\sum^{K-1}_{k=0} \big| \log Df (f^k (g^N (z))) - \log Df (f^k(z)) \big| \leq \var (\log Df) \leq V.$$
We thus conclude that \, $\log (D (f^K g^N f^{-K}) (z_K))  > V,$ \, as announced.


\section{Distortion in class $C^{1+\alpha}$ for $\alpha < 1/2$} 

In this section, we start by briefly recalling the construction of the group $\Gamma$ with a distorted element $\bar{f}$ considered in \cite{Na}. 
Next, we proceed to smooth the action of $\Gamma$ in order to achieve any differentiability class $C^{1+\alpha}$ for $\alpha < 1/2$. Upgrading 
$\alpha$ to any number less than 1 will require the introduction of an extra element plus a tricky new computation, and will be carried out in 
the next section. 

Start with the vector fields $\hat{\mathcal{X}}$ and $\mathcal{X}$ on the real line 
whose time-1 maps are, respectively, 
$$\hat{F} 
:= \hat{\mathcal{X}}^1
: x \mapsto 2x \qquad { \mathrm{and} } \qquad F
:= \mathcal{X}^1
: x \mapsto x + 1.$$
Let $\varphi \! : \mathbb{R} \to \,  ]0,1[$  be a $C^{\infty}$ diffeomorphism such that  \, 
$\hat{\mathcal{Y}} := \varphi_* (\hat{\mathcal{X}})$ \, and \, $\mathcal{Y} := \varphi_* (\mathcal{X})$ \, 
extend to the endpoints of $[0,1]$ as infinitely flat vector fields. Denote $\hat{f} := \hat{\mathcal{Y}}^1$ and 
$f := \mathcal{Y}^1$, which we view as diffeomorphisms of $[-1,2]$ that coincide with the identity outside $[0,1]$. 
The affine relation $\hat{f} f \hat{f}^{-1} = f^2$ yields that $\| f^n \| = O (\log (n))$; in particular, 
$f$ has vanishing asymptotic distortion.

Let $x_0 := \varphi (0)$ and, for each $k \!\in\! \mathbb{Z}$, let $x_k := f^k (x_0) = \varphi (k)$. Denote 
also $x_{-1/2} := \varphi (-1/2)$ and $x_{-3/4} := \varphi (-3/4)$. Let $\varphi_0$ the affine diffeomorphism 
sending $I := [x_0,x_1]$ onto $[0,1]$, and let $g := \varphi_0^{-1} f \varphi_0$. This can be extended to $[-1,2]$ 
by the identity outside $I$. 

We next define two diffeomorphisms $\hat{h}$ and $h$ as follows:

\vspace{0.18cm}

\noindent (i) They act by the identity outside $[0,1]$.

\vspace{0.18cm}

\noindent (ii) On each interval $I_k := f^k (I)$, the diffeomorphism $\hat{h}$ (resp. $h$) coincides with the 
$s_k$-time map (resp. $t_k$-time map) of the flow of the vector field $f^k_* \big(\varphi_0^* (\hat{\mathcal{Y}} ) \big)$ 
(resp. $f^k_* \big(\varphi_0^* (\mathcal{Y})) $\big). 

\vspace{0.18cm}

\noindent Here, $s_k$ and $t_k$ are sequences of real numbers such that:

\vspace{0.18cm}

\noindent (iii) If $2^{i-1} \leq k < 2^{i}$ for a certain positive even integer $i$, then 
$$s_k := \log_2 \left( 1 - \frac{1}{\sqrt{\ell_{i/2}}} \right) \qquad \mbox{ and } \qquad t_k := \frac{1}{ \sqrt{\ell_{i/2}} },$$
where $\ell_j$ is a prescribed sequence of positive integers diverging to infinity to be fixed below.

\vspace{0.18cm} 

\noindent (iv) Otherwise, $s_k = t_k := 0$.

\vspace{0.18cm}

Finally, 
we let $\psi$ be a $C^{\infty}$ diffeomorphism of $[-1,2]$ such that:

\vspace{0.18cm}

\noindent (v) $\psi$ coincides with the identity on $[x_{-1/2} , x_0],$


\vspace{0.18cm}

\noindent (vi) $\psi (x_{-3/4}) = 0$ and $\psi (x_1) = 1$. 

\vspace{0.18cm}

\noindent The group we consider is  $\Gamma := \langle \hat{f}, f, g, \hat{h}, h, \psi \rangle$. The 
computations of \cite{Na} show the following relation for certain powers of $\bar{f} := fg$ (which justifies the construction): 
\begin{equation}\label{eq-identity}
(\bar{f})^{2^{i-1}} =  f^{1+n/2} \hat{f}^i [\hat{f}^{-i}f^{-n}hf^n\hat{f}^i ,  \psi \hat{f}^{-i}f^{-n}\hat{h}f^n\hat{f}^i\psi^{-1} ]^{\ell_{i/2}} \hat{f}^{-i}  f^{-1} 
\end{equation}
where $i$ is an even (positive) integer and $n:= 2^i$. Roughly, this works as follows: set 
$$a_n:=\hat{f}^{-i}f^{-n}hf^n\hat{f}^i \quad \mbox{ and } \quad b_n:= \psi \hat{f}^{-i}f^{-n}\hat{h}f^n\hat{f}^i\psi^{-1}.$$ 
One easily checks that 
$$\text{supp} (a_n) \subset[0,x_{-3/2}]\cup[x_{-1/2},x_0]\cup[x_1,1] 
\quad \mbox{ and } \quad \text{supp} (b_n) \subset[-1,0]\cup[x_{-1/2},x_0]\cup[1,2].$$
Thus, the commutator \, $c_n:=[a_n,b_n] = a_n b_n a_n^{-1} b_n^{-1}$ \, is supported on $[x_{-1/2},x_0]$, hence 
the conjugate $\hat{f}^{i}c_n\hat{f}^{-i}$ is supported on $ [x_{-2^{i-1}},x_0] $. Besides, on each $[x_k,x_{k+1}]\subset [x_{-2^{i-1}},x_0]$, 
this conjugate $\hat{f}^{i}c_n\hat{f}^{-i}$ coincides with the time-$\frac{1}{\ell_{i/2}}$ map of the flow of $f^k_\ast(\varphi_0^\ast(\mathcal{Y}))$.
Moreover, the restriction of the map 
$$h_{n/2} := (f^{-n/2}gf^{n/2}) \cdots (f^{-2} g f^2) (f^{-1}gf) = f^{-n/2} (f^{-1} \bar{f}^{n/2} f)$$ 
to each $[x_k,x_{k+1}]\subset [x_{-2^{i-1}},x_0]$ equals the time-1 map of the flow of $f^k_\ast(\varphi_0^\ast(\mathcal{Y}))$. This implies that 
$$h_{n/2} = (\hat{f}^i c_n \hat{f}^{-i})^{\ell_{i / 2}},$$ 
which corresponds to (\ref{eq-identity}). 

\vspace{0.2cm}

Since $\| f^n \| = O (\log(n))$, identity (\ref{eq-identity}) implies that $\| (\bar{f})^{2^{i-1}}\| = O (i \, \ell_{i/2})$. 
Therefore, $\bar{f}$ is distorted provided $\ell_{j}$ grows to infinite in such a way that 
\begin{equation}\label{eq-ell}
\lim_{n \to \infty} \frac{\log (n) \, \ell_i}{n}  
= \lim_{i \to \infty} \frac{i \, \ell_i }{2^i} = 0
\end{equation}

The maps $\hat{f}, f, g, \psi$ above are obviously smooth. However, regularity for $\hat{h}, h$ is more subtler. Indeed, 
their $C^1$ smoothness is ensured by the conditions $s_n \to 0$ and $t_n \to 0$ as $|n| \to \infty$ (which are equivalent 
to $\ell_j \to \infty$ as $j \to \infty$) together with the next lemma. This is strongly inspired on the work of Kopell (a proof 
together with a discussion appears in the Appendix of \cite{Na}). Below, by $C^{1+Lip}$ we refer to maps with Lipschitz derivative.

\vspace{0.1cm}

\begin{lem}\label{lem-gen-kopell}
Let $f \in \Diff^{1+Lip}_+([0,1])$ be such that $f(x) \neq x$ for all $x \in (0,1)$. Fix $x_0\in (0,1)$, and let $x_n := f^n(x_0)$. 
Let $(g_n)_{n \in \mathbb{Z}}$ be a sequence of $C^1$ diffeomorphisms of the interval of endpoints $x_0,x_1$ such that 
$Dg_n(x_0) = Dg_n(x_1)=1$ for all $n \in \mathbb{Z}$. Let $g: (0,1)\rightarrow (0,1)$ be the diffeomorphism whose 
restriction to each interval of endpoints $x_n , x_{n+1}$ coincides with $f^ng_nf^{-n}$. If  $g_n \rightarrow Id$ \, in $C^1$ 
topology as $n \to \infty$, then $g$ extends to a $C^1$ diffeomorphism of $[0,1]$ by letting $f(0) = 0$ and $f(1)=1$.
\end{lem}

\vspace{0.1cm}

Unfortunately, as we will see in the Appendix, this fails to extend to class $C^{1+\alpha}$. Because of this, we need to go 
into more explicit computations for our example. Although these are difficult to handle, the following key elementary lemma 
taken from the work of Pixton \cite{pixton} and Tsuboi \cite{tsuboi} will be enough for us.

\vspace{0.1cm}

\begin{lem} \label{lema-a-la-Tsuboi}
{\em Given a $C^2$ vector field $\mathcal{X}$ on an interval $[0,a]$, denote 
\, $C_1 := \| D \mathcal{X} \|$ \, and \, $C_2 := \| D^2 \mathcal{X} \| $. \,
If $f^t$ denotes its flow, then, for all $t \geq 0$,}
$$\big\| D \log Df^t \big\| \leq \frac{C_2}{C_1} \, \big( e^{C_1 t} -1 \big).$$
\end{lem}

\noindent{\bf Proof.} Taking derivatives on the equality $d f^t / dt = \mathcal{X} \circ f^t$, we deduce
$$\frac{d}{dt} D f^t = D \mathcal{X} (f^t) \cdot Df^t,$$
hence
$$\frac{d}{dt} \log Df^t = \frac{\frac{d}{dt} Df^t}{Df^t} = D \mathcal{X} (f^t).$$
Since $f^0 = Id$, we conclude
$$ \log Df^t  = \int_0^t D \mathcal{X} (f^s) \, ds.$$
Since \, $| D \mathcal{X} | \leq C_1$, \, this yields \, $| \log Df^t | \leq C_1 t$, \, hence \, $|D f^t| \leq e^{C_1 t}$. \, 
Moreover,
$$D \log Df^t = D \left( \int_0^t D \mathcal{X} (f^s) \, ds \right) = \int_0^t D^2 \mathcal{X} (f^s) \cdot Df^s \, ds.$$
Since \, $| D^2 \mathcal{X} | \leq C_2$, \, we conclude
$$| D \log Df^t | \leq C_2 \int_0^t |Df^s| \, ds \leq C_2 \int_0^t e^{C_1 s} ds = \frac{C_2}{C_1} \, \big( e^{C_1 t} -1 \big),$$
as announced. 
$\hfill\square$

\vspace{0.35cm}

We now turn into very long computations that will allow us to ensure that the resulting maps $h, \hat{h}$ built via 
the procedure above are $C^{1+\alpha}$ diffeomorphisms for well chosen $\varphi$ and $\ell_i$ that respect all the 
properties we have imposed. This will close the proof of our theorem.

In order to simplify these computations, let us remind the {\em chain rules} for different 
derivatives of maps between 1-dimensional spaces, namely logarithmic ($L$), affine ($A$), and Schwarzian ($S$):  
$$L (f) := \log Df, \quad A(f) := DL(f) = \frac{D^2 f}{D f} , \quad S (f) := DA(f) - \frac{A(f)^2}{2} = \frac{D^3 f}{D f}  - \frac{3}{2} \left( \frac{D^2 f}{D f} \right)^2.$$
These are listed below: 
\begin{equation}\label{log}
L(fg) = L(g) + L(f) \circ g,
\end{equation}
\begin{equation}\label{afin}
A (fg) = A (g) + A(f) \circ g \cdot Dg,
\end{equation}
\begin{equation}\label{sw}
S (fg)= S (g) + S(f) \circ g \cdot (Dg)^2 . 
\end{equation}

We let $\varphi \! : (0,1) \to \mathbb{R}$ be a $C^{\infty}$ diffeomorphism such that, 
for a small-enough $\delta > 0$ 
$$\varphi(x)= \left\{ \begin{array}{lcc}
            - \text{exp}{(\text{exp}{(1/x)})} &   \mbox{if}  & 0<x\leqslant \delta 
             \\  \text{exp}{(\text{exp}{(1/(1-x))})} &  \mbox{if}  & 1-\delta \leqslant x <1
             \end{array}
   \right.$$    
If we denote  $\mathcal{Z} := \varphi_0^* (\hat{\mathcal{Y}} )$, then we need to control $f^n_* (\mathcal{Z})$. 

\vspace{0.32cm}

\noindent{\bf An estimate for lengths of fundamental domains.} Let us come back to the group 
$\Gamma = \langle \hat{f}, f, g, \hat{h}, h, \psi \rangle$. 
Recall that $I$ denotes the interval $[x_0,x_1]$. We claim that, for a certain constant $C > 0$,  
\begin{equation}\label{eq-long}
 |f^n (I)|  = O \left( \frac{C}{n \log(n) \, (\log (\log (n)))^2} \right).
\end{equation}
This is checked via a direct computation. Namely, for a large-enough $n$, 
$$ | f^n (I)| = \varphi^{-1} (n+1) - \varphi^{-1} (n) = \frac{1}{\log \log (n)} - \frac{1}{\log \log (n+1)}. $$
Since the derivative of \, $x \mapsto 1/\log\log(x)$ \, is \, $1 / [x \, \log (x) \, (\log\log(x))^2]$, \, a direct 
application of the Mean Value Theorem yields the desired estimate (\ref{eq-long}).

\vspace{0.32cm}

\noindent{\bf Estimates for the vector field and its derivative.} Notice that 
$$f^n_* (\mathcal{Z}) (x) = [Df^n(f^{-n}(x))] \, \mathcal{Z}(f^{-n}(x)), \hspace{0,5cm} x\in I_n := f^n (I) .$$
Taking derivatives, we obtain
$$D (f^n_\ast(\mathcal{Z}))(x) = \frac{D^2f^n(f^{-n}(x))}{Df^{n}(f^{-n}(x))} \, \mathcal{Z}(f^{-n}(x)) + D \mathcal{Z}(f^{-n}(x)).$$ 
Now, using the chain rule (\ref{afin}), this yields
$$ D (f^n_\ast(\mathcal{Z}))(x) 
= \sum_{i=0}^{n-1}\left( \frac{D^2f(f^{i-n}(x))}{Df(f^{i-n}(x))}\right)Df^{i}(f^{-n}(x)) \, \mathcal{Z}(f^{-n}(x)) + D \mathcal{Z}(f^{-n}(x)).$$ 
Thus,  letting 
$$C' := \left\| \frac{D^2 f}{Df} \right\| \| \mathcal{Z} \|, \qquad C'' := \| D \mathcal{Z} \|, $$
we obtain
$$ | D (f^n_\ast (\mathcal{Z}))(x) | \leq C' \sum_{i=0}^{n-1}  Df^{i}(f^{-n}(x))  + C''.$$
We claim that the sum above is uniformly bounded (independently of $n$ and $x \in I_n$), so that 
\begin{equation}\label{eq-est-der}
| D (f^n_\ast (\mathcal{Z}))(x) | \leq C
\end{equation}
for a certain constant $C$. Indeed, a standard control of distortion argument yields that $Df^{i}(f^{-n}(x))$ is of the order of 
$|I_i| / |I_0|$, hence 
$$ \sum_{i=0}^{n-1}  Df^{i}(f^{-n}(x)) \sim  \sum_{i=0}^{n-1}  |I_i| \leq 1.$$

\vspace{0.32cm}

\noindent{\bf Estimates for the second derivative.} We now claim that, for a certain constant $C>0$ and all $x \in I_n$,
\begin{equation}\label{est-D2}
\big| D^2  (f^n_\ast(\mathcal{Z})) (x) \big| \leq C \, n \log (n) \, (\log \log (n))^2.
\end{equation}
To show this, notice that, from
$$D (f^n_\ast(\mathcal{Z}))(x) = A(f^n) (f^{-n}(x)) \, \mathcal{Z}(f^{-n}(x)) + D \mathcal{Z}(f^{-n}(x)),$$ 
we obtain
$$D^2 (f^n_\ast(\mathcal{Z}))(x) 
= \left[ DA(f^n) \cdot \mathcal{Z}  + A(f^n)  \, D \mathcal{Z} + D^2 \mathcal{Z} \right] \circ f^{-n} (x) \cdot Df^{-n}(x).$$ 
which is equal to 
$$\big[ \big( S(f^n) + \frac{1}{2} A(f^n)^2\big) \cdot \mathcal{Z}  + A(f^n)  \, D \mathcal{Z} + D^2 \mathcal{Z} \big] \circ f^{-n} (x) \cdot Df^{-n}(x).$$ 

Let us analise each term entering in this expression. First, 
by (\ref{eq-long}) and the control of distortion argument above, 
$$Df^{-n}(x) = 1 / D f^n (f^{-n}(x)) = O \big( n \log (n) \, (\log \log (n))^2 \big).$$

We next claim that $A(f^n)$ is uniformly bounded on $I_0$. Indeed, letting $C:= \| A (f) \|$, the chain rule (\ref{afin}) yields  
$$A (f^n) = \sum_{i=0}^{n-1} A (f) \circ f^i \cdot Df^i \leq C \sum_{i=0}^{n-1} Df^i.$$
The control of distortion argument above shows that the last sum is bounded from above by a constant, hence the claim.

Since $\mathcal{Z}$, $D \mathcal{Z}$ and $D^2 \mathcal{Z}$ are obviously uniformly bounded, to show (\ref{est-D2}) 
it remains to check that $S (f^n) (f^{-n} (x))$ is uniformly bounded. To see this, we use the chain rule (\ref{sw}):
$$Sf^n(f^{-n}(x)) = \sum_{i=0}^{n-1}Sf(f^{i-n}(x))(Df^{i}(f^{-n}(x)))^2.$$ 
This implies
$$ \big| Sf^n(f^{-n}(x)) \big| \leq C \, \sum_{i=0}^{n-1} \big( Df^i (f^{-n}(x)) \big)^2,$$
and the last sum can be estimated as it was done before. (The sum here is even smaller since it involves the squares of the derivatives.)

\vspace{0.32cm}

\noindent{\bf Estimates for the maps.} We are now in position to check that the group $\Gamma$ is made of $C^{1+\alpha}$ diffeomorphisms 
for $\alpha < 1/2$ and $\ell_i$ of order $n / \log(n)^2$. (Notice that, according to (\ref{eq-ell}), the element $\bar{f}$ is  distorted in $\Gamma$ 
for this choice.) Notice that this is obvious for all the generators except $h$ and $\hat{h}$. The estimates for these two elements are similar, 
so that we only deal with $h$. Besides, we may deal with $\log Dh$ instead of $Dh$, since the condition ``$Dh$ is of class $C^{\alpha}$'' is 
equivalent to that ``$\log Dh$ is of class $C^{\alpha}$''.

We need to check that there exists a uniform bound $B$ for expressions of type
$$\frac{|\log Dh(y) - \log Dh(x)|}{|y-x|^{\alpha}}$$
for all points $x< y$ in the same interval $I_n$. Indeed, having such an estimate, one can easily 
treat the case of arbitrary pairs $x < y$ just  by noticing that at each endpoint of an interval of the form above, 
the derivative of $h$ equals 1. Namely, letting $z_1$ (resp. $z_2$) be such an endpoint that is immediately 
to the right of $x$ (resp. to the left of $y$), one has
\begin{eqnarray*}
\frac{|Dh(y) - Dh(x)|}{|y-x|^{\alpha}} 
&\leq& \frac{|Dh(y) - Dh(z)|}{|y-x|^{\alpha}} + \frac{|Dh(z) - Dh(x)|}{|y-x|^{\alpha}} \\
&\leq& \frac{|Dh(y) - Dh(z)|}{|y-z|^{\alpha}} + \frac{|Dh(z) - Dh(x)|}{|z-x|^{\alpha}} \\
&\leq& 2 \, B.
\end{eqnarray*}

Now, for all $ z \in I_n$ (with $n \geq 0$), Lemma \ref{lema-a-la-Tsuboi} and estimate (\ref{eq-est-der}) yield, for $t_n$ small enough, 
$$D (\log Dh) (z) \leq 2 \, \| D^2 f^n_* (\mathcal{Z}) \| \, t_n.$$ 
By estimate (\ref{est-D2}), this implies, for a certain constant $C>0$,  
\begin{equation}\label{eq-extra}
D (\log Dh) (z) \leq 2 \, C \, n \log (n) \, (\log(\log(n)))^2 \, t_n.
\end{equation}
Moreover, for $x,y$ in $I_n$, 
$$\frac{ | \log Dh (y) - \log Dh (x)| }{|y-x|^{\alpha}} = \frac{ | \log Dh (y) - \log Dh (x)| }{|y-x|} \, |y-x|^{1-\alpha} = D (\log Dh (z)) \, |y-x|^{1-\alpha}$$
for a certain point $z \in I_n$. By (\ref{eq-extra}), this yields
\begin{equation}\label{eq-z}
\frac{ | \log Dh (y) - \log Dh (x)| }{|y-x|^{\alpha}}  
\leq 
2 \, C \, n \log (n) \, (\log(\log(n)))^2 \, t_n \, \left[  \frac{C}{ n \log (n) (\log (\log (n)))^2}, \right] ^{1-\alpha}.
\end{equation}
Since \, $t_n = 1/ \sqrt{\ell_{i/2}} \leq C \, \log(n) / \sqrt{n}$, \, we finally obtain
$$\frac{ | \log Dh (y) - \log Dh (x)| }{|y-x|^{\alpha}} 
\leq 2 \, C' \, n \log (n) \, (\log(\log(n)))^2 \, \frac{\log(n)}{n^{1/2}} \, \frac{1}{[n \log (n) (\log(\log(n)))^2]^{1-\alpha}}.$$
To get the desired upper bound $B$, it suffices that the total exponent of $n$ in the expression above is negative. 
Since this exponent equals  \, $1 - 1/2 - (1-\alpha) = \alpha - 1/2$, \, this condition reduces to $\alpha < 1/2$, which is our hypothesis. 


\section{Distortion in class $C^{1+\alpha}$ for $\alpha < 1$} 

It is unclear whether the previous action can be smoothed beyond the class $C^{3/2}$ (compare \cite{int4,int1,int2,int3}). To achieve a larger 
differentiability class, we will need to accelerate the distorted 
behavior of $\bar{f}$, which will allow us to consider smaller integration times for the flows of vector fields (in concrete terms, we will increase 
the sequence $\ell_i$). This will be crucial to improve the regularity from $\alpha < 1/2$ to any $\alpha < 1$.

\vspace{0.32cm}

\noindent{\bf Adding an extra element.} We consider the 
map $\tilde{h}$ acts by the identity outside the intervals $I_k$, and that on each such interval coincides with the $r_k$-time of the 
time flow of the vector field $f^k_* (\varphi_0^* (\hat{\mathcal{Y}}))$, where $r_k := 1 / \sqrt{\ell_{i/2}}$ for $2^{i-1} \leq k < 2^{i}$ and 
$r_k := 0$ otherwise. Notice that $\tilde{h}$ is very similar to $\hat{h}$. (Actually, we could perform the computations that 
follow using $\hat{h}$ instead of $\tilde{h}$, but this would become much harder.)

Then we let \, $d_n := \hat{f}^{-i}f^{-n}\tilde{h}f^n\hat{f}^{i}$ \, for $n = 2^i$, where $i$ is an even 
integer. We have \, $\mathrm{supp} (d_n) \subset [0,x_{-3/2}]\cup[x_{-1/2},x_0]\cup[x_1,1]$. \, 
Since \, $\mathrm{supp} (c_n) \subset [x_{-1/2},x_0]$, \, for every integer $L_i \geq 1$, 
the support of \, $d_n^{L_i} c_n d_n^{-L_i}$ \, is also 
contained in $[x_{-1/2}, x_0]$, thus the support of 
$$\hat{f}^{-i} d_n^{L_i} c_n d_n^{-L_i} \hat{f}^i 
= (\hat{f}^{-i} d_n^{L_i} \hat{f}^{i}) (\hat{f}^{-i} c_n \hat{f}^{i}) (\hat{f}^{-i} d_n^{-L_i} \hat{f}^{i})$$
is contained in $[x_{-2^{i-1}},x_0]$.

Now recall that, on each $[x_k,x_{k+1}]\subset [x_{-2^{i-1}},x_0]$, the conjugate $\hat{f}^{i}c_n\hat{f}^{-i}$ coincides with the 
time-$\frac{1}{\ell_{i/2}}$ map of the flow of $f^k_\ast(\varphi_0^\ast(\mathcal{Y}))$. Moreover, by construction, on the same interval, the  
conjugate $\hat{f}^{i} d_n\hat{f}^{-i}$ coincides with the time-$\frac{1}{\sqrt{\ell_{i/2}}}$ map of the flow of $f^k_\ast(\varphi_0^\ast(\hat{\mathcal{Y}}))$. 
By the affine relation, still on the same interval, the map $\hat{f}^{-i} d_n^{L_i} c_n d_n^{-L_i} \hat{f}^i$ lies in the flow of 
$f^k_\ast(\varphi_0^\ast(\mathcal{Y}))$, and arises at time
$$\frac{2^{\frac{L_i}{\sqrt{\ell_{i/2}}}}}{\ell_{i/2}}.$$
If $L_i := \sqrt{\ell_{i/2}} \, \log_2 (\ell_{i/2})$ (which will be chosen to be an integer number), 
then this quantity equals 1. Therefore, for this choice,  $\hat{f}^{-i} d_n^{L_i} c_n d_n^{-L_i} \hat{f}^i$ coincides with $h_{n/2}$.  

\vspace{0.32cm}

\noindent{\bf The distortion estimate.} The identity \, $h_{n/2} = \hat{f}^{-i} d_n^{L_i} c_n d_n^{-L_i} \hat{f}^i$ \,  
implies that, in the new group 
$\tilde{\Gamma} := \langle \hat{f}, f, g, \hat{h}, h, \tilde{h}, \psi \rangle$, we have the estimate 
$$\lVert h_{n/2} \rVert \leqslant 2\lVert \hat{f}^i\rVert + 2 \, L_i \lVert d_n\rVert+\lVert c_n\rVert\leqslant 2i + 2 \, L_i \, 
(2i+1+2\lVert f^n\rVert) + 8 \, (1+i+\lVert f^n\rVert).$$ 
Since $\| f^n \| = O (\log (n)) = O (i)$, we conclude that 
$$\lVert h_{n/2} \rVert = O \left( i \, \sqrt{\ell_{i/2}} \, \log (\ell_{i/2}) \right).$$
Since $h_{n/2} = f^{-n/2} (f^{-1} \bar{f}^{n/2} f) $ and $\| f^{n/2} \| = O (\log(n))$, this yields
$$\| \bar{f}^{n/2} \| = 2 + \| f^{n/2} \| +  \| h_{n/2} \| = O \left( i \, \sqrt{\ell_{i/2}} \, \log (\ell_{i/2}) \right).$$

Notice that the last estimate is much better than what we had in the group $\Gamma$ of the previous section. 
In there, $\| \bar{f}^{n/2}\|$ 
was of the order $O (i \, \ell_{i/2})$, hence, $\bar{f}$ was distorted provided the growth of $\ell_j$ was smaller than exponential. In the new setting, 
that is, in the modified group $\tilde{\Gamma}$, the diffeomorphism $\bar{f}$ is distorted whenever the condition below is satisfied (recall that $n = 2^i$):
\begin{equation}\label{eq-new-dist}
\lim_{n \to \infty}  \frac{i \, \sqrt{\ell_{i/2}} \, \log (\ell_{i/2})}{2^i} = 0.
\end{equation}
 
\vspace{0.32cm}

\noindent{\bf Checking regularity.} We thus choose a new sequence $\ell_i$ so that condition (\ref{eq-new-dist}) holds and 
$\sqrt{\ell_{i/2}} \, \log_2 (\ell_{i/2})$ is an integer number. This can be achieved for a sequence of type 
$$\sqrt{\ell_{i/2}}\sim \frac{n}{\log (n)^{3}},$$
that we fix from now on. 
With such a choice, we claim that $\tilde{\Gamma}$ is a group of 
$C^{1+\alpha}$ diffeomorphisms. Again, this is obvious for all generators except $h,\hat{h},\tilde{h}$, and for these three 
elements the computations are the exact same, because each of the sequences $r_n, s_n, t_n$ is equivalent to $1/ \sqrt{\ell_{i/2}}$.  
We thus write everything only for $h$. Remind estimate (\ref{eq-z}):
$$\frac{ | \log Dh (y) - \log Dh (x)| }{|y-x|^{\alpha}}  
\leq 
2 \, C \, n \log (n) \, (\log(\log(n)))^2 \, t_n \, \left[  \frac{C}{ n \log (n) (\log (\log (n)))^2}, \right] ^{1-\alpha}.$$
With the new estimate for $t_n$, this becomes
$$\frac{ | \log Dh (y) - \log Dh (x)| }{|y-x|^{\alpha}}  
\leq 
2 \, C \, n \log (n) \, (\log(\log(n)))^2 \, \frac{\log (n)^{3}}{n} \, \left[  \frac{C}{ n \log (n) (\log (\log (n)))^2}, \right] ^{1-\alpha}.$$
The expression on the right is of order
$$O \left( \frac{\log(n)^{3+\alpha} \, (\log \log(n))^{2\alpha}}{n^{1-\alpha}} \right),$$
which converges to $0$ as $n$ goes to infinite. This allows showing that $h$ is a $C^{1+\alpha}$ diffeomorphism as it was 
done in the previous section.
 

\section{Appendix: no strong Kopell's lemma in class $C^{1+\alpha}$}

The goal of this Appendix is to give an example of the phenomenon announced just after Lemma \ref{lem-gen-kopell}.
Namely, there exist: 

\vspace{0.1cm}

\noindent -- a $C^{\infty}$ diffeomorphism of $[0,1]$ such that $f(x) > x$ for all $x \in (0,1)$ with  a fundamental domain $[x_0,x_1]$ (where $x_1 := f(x_0)$), 

\vspace{0.1cm}

\noindent and 

\vspace{0.1cm}

\noindent -- a sequence $(g_n)_{n \in \mathbb{Z}}$ of $C^{1+\alpha}$ diffeomorphisms of $[x_0,x_1]$, 

\vspace{0.1cm}

\noindent such that \, 
$Dg_n (x_0) = Dg_n (x_1) = 1$ \, for all  $n$, one has the convergence $g_n \to Id$ in $C^{\infty}$ topology, 
but the $C^1$ diffeomorphism $g \! : [0,1]\rightarrow [0,1]$ defined as 
\begin{equation}\label{def-g_n}
g|_{f^n([x_0,x_1])}: =    f^n g_nf^{-n}|_{f^n([x_0,x_1])} 
\end{equation}
is not $C^{1+\alpha}$ for any $\alpha \in (0,1)$.

\vspace{0.1cm}

Our diffeomorphism $f$ is 
$$f(x) = \frac{2x}{x+1}.$$
Notice that
$$f^n(x)=\frac{2^nx}{(2^{n}-1)x+1}.$$ 
Hence, for all $x,y$,
$$\big| f^n (x) - f^n (y) \big| = \left| \frac{2^n x}{(2^{n}-1)x+1} - \frac{2^n y}{(2^{n}-1)y + 1} \right| 
=  \frac{2^n |x-y|}{[(2^{n}-1)x+1] \cdot [(2^{n}-1)y+1]},$$
which yields
\begin{equation}\label{eq-distancia}
\big| f^n (x) - f^n (y) \big| \leq \frac{C}{2^n}
\end{equation}
for a certain universal constant $C$ provided $x,y$ both belong to a compact subinterval of 
$(0,1)$ (say $[1/2,2/3]$). 

Now consider the points $x_0 := 1/2$ and $x_1 := f (x_0) = 2/3$. Notice that $x_1-x_0 = 1/6$. Then let
$$a:= \frac{1}{2} + \frac{1}{30} = \frac{8}{15}, \quad b :=  \frac{1}{2} + \frac{2}{30} = \frac{17}{30}, \quad 
b' :=  \frac{1}{2} + \frac{3}{30} = \frac{3}{5}, \quad b'' := \frac{1}{2} + \frac{4}{30} = \frac{19}{30}.$$
Let $\varrho \! : [a,2/3] \to [-1/2,1/2]$ be a $C^{\infty}$ function such that 
$$\varrho (a) = \varrho (b') = \varrho (2/3) = 0, \qquad 
\varrho (b) = \frac{1}{2}, \qquad \varrho (b'') = -\frac{1}{2}.$$
Assume also that $\varrho$ is strictly increasing on $[a,b]$ and $[b'',3/2]$, strictly decreasing on $[b,b'']$, 
infinitely flat at $a$ and $3/2$, and its graph is symmetric with respect to the point $(b',0)$. 
Let $\rho_n \! : [1/2,2/3]$ be the function that is identically equal  to $1$ on $[1/2,a]$ 
and whose restriction to $[a,2/3]$ coincides with $1 + \varrho / n$. 

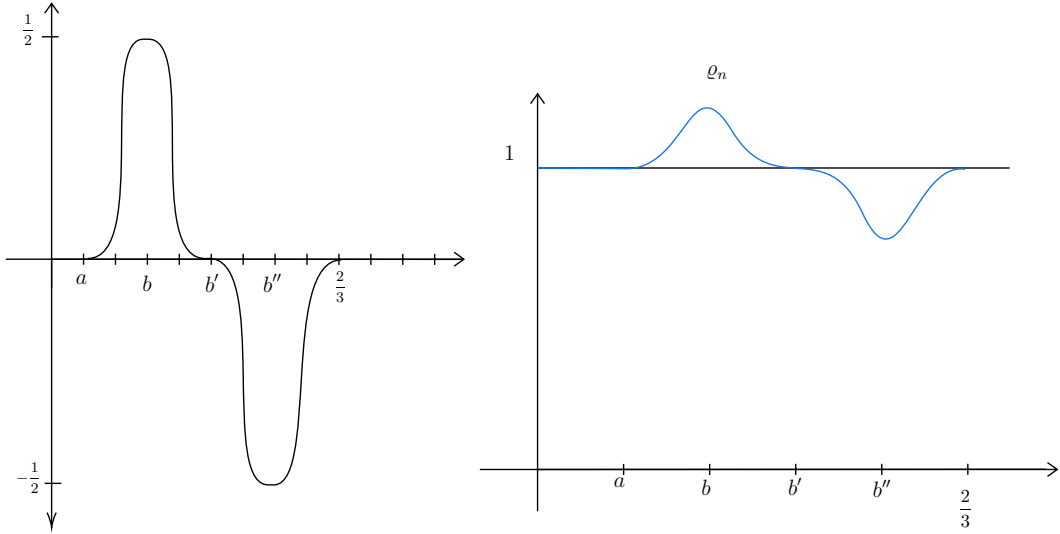
\begin{center}

\tikzset{every picture/.style={line width=0.75pt}} 

\scalebox{0.7}{\begin{tikzpicture}[x=0.75pt,y=0.75pt,yscale=-1,xscale=1]

\draw  (50,248.38) -- (380.02,248.38)(83,64) -- (83,268.86) (373.02,243.38) -- (380.02,248.38) -- (373.02,253.38) (78,71) -- (83,64) -- (88,71)  ;
\draw    (109.02,247.86) .. controls (154.02,247.86) and (115.02,83.86) .. (151.02,89.86) ;
\draw    (195.02,247.86) .. controls (148.02,251.23) and (190.02,84.23) .. (151.02,89.86) ;
\draw    (291.91,248.85) .. controls (246.91,248.79) and (276.99,417.05) .. (241,411) ;
\draw    (195.02,247.86) .. controls (242.02,244.57) and (201.99,416.58) .. (241,411) ;
\draw    (83,248.38) -- (83,440) ;
\draw [shift={(83,442)}, rotate = 270] [color={rgb, 255:red, 0; green, 0; blue, 0 }  ][line width=0.75]    (10.93,-3.29) .. controls (6.95,-1.4) and (3.31,-0.3) .. (0,0) .. controls (3.31,0.3) and (6.95,1.4) .. (10.93,3.29)   ;
\draw    (83,248.38) -- (334.5,248.48) -- (346.5,248.49) -- (377.5,248.5) (106,244.39) -- (106,252.39)(129,244.4) -- (129,252.4)(152,244.41) -- (152,252.41)(175,244.42) -- (175,252.42)(198,244.43) -- (198,252.43)(221,244.43) -- (221,252.43)(244,244.44) -- (244,252.44)(267,244.45) -- (267,252.45)(290,244.46) -- (290,252.46)(313,244.47) -- (313,252.47)(336,244.48) -- (336,252.48)(359,244.49) -- (359,252.49) ;
\draw    (76,88) -- (88,88) ;
\draw    (78,410) -- (90,410) ;

\draw (60,70.4) node [anchor=north west][inner sep=0.75pt]  [font=\scriptsize] [align=left] {$\displaystyle \frac{1}{2}$};
\draw (56,393.4) node [anchor=north west][inner sep=0.75pt]  [font=\scriptsize] [align=left] {$\displaystyle -\frac{1}{2}$};
\draw (99,257.86) node [anchor=north west][inner sep=0.75pt]   [align=left] {$\displaystyle a$};
\draw (147,258.86) node [anchor=north west][inner sep=0.75pt]  [font=\normalsize] [align=left] {$\displaystyle b$};
\draw (192.02,257.86) node [anchor=north west][inner sep=0.75pt]   [align=left] {$\displaystyle b'$};
\draw (234.02,257.86) node [anchor=north west][inner sep=0.75pt]   [align=left] {$\displaystyle b''$};
\draw (285,255.86) node [anchor=north west][inner sep=0.75pt]  [font=\scriptsize] [align=left] {$\displaystyle \frac{2}{3}$};

\end{tikzpicture}

\tikzset{every picture/.style={line width=0.75pt}} 

\begin{tikzpicture}[x=0.75pt,y=0.75pt,yscale=-1,xscale=1]

\draw  (65.02,292.71) -- (481,292.71)(106.62,21.89) -- (106.62,322.8) (474,287.71) -- (481,292.71) -- (474,297.71) (101.62,28.89) -- (106.62,21.89) -- (111.62,28.89)  ;
\draw    (106.62,292.71) -- (472.95,292.57) (168.61,288.68) -- (168.62,296.68)(230.61,288.66) -- (230.62,296.66)(292.61,288.64) -- (292.62,296.64)(354.61,288.61) -- (354.62,296.61)(416.61,288.59) -- (416.62,296.59) ;
\draw [color={rgb, 255:red, 0; green, 0; blue, 0 }  ,draw opacity=1 ]   (106.59,75.32) -- (446.76,75.32) ;
\draw [color={rgb, 255:red, 19; green, 118; blue, 233 }  ,draw opacity=1 ]   (106.59,75.32) -- (167.65,75.78) ;
\draw [color={rgb, 255:red, 19; green, 118; blue, 233 }  ,draw opacity=1 ]   (167.65,75.78) .. controls (214.61,78.98) and (215.96,-2.08) .. (246.82,48.75) .. controls (277.68,99.59) and (313.91,50.13) .. (340.75,107.84) .. controls (367.59,165.55) and (378.32,70.74) .. (414.55,76.23) ;

\draw (80.98,57.82) node [anchor=north west][inner sep=0.75pt]   [align=left] {$\displaystyle 1$};
\draw (160.15,296.9) node [anchor=north west][inner sep=0.75pt]   [align=left] {$\displaystyle a$};
\draw (223.22,298.28) node [anchor=north west][inner sep=0.75pt]   [align=left] {$\displaystyle b$};
\draw (285.63,296.9) node [anchor=north west][inner sep=0.75pt]   [align=left] {$\displaystyle b'$};
\draw (346.7,298.28) node [anchor=north west][inner sep=0.75pt]   [align=left] {$\displaystyle b''$};
\draw (408.57,305.83) node [anchor=north west][inner sep=0.75pt]  [font=\footnotesize] [align=left] {$\displaystyle \frac{2}{3}$};
\draw (226.61,0.11) node [anchor=north west][inner sep=0.75pt]   [align=left] {$\displaystyle \varrho _{n} \ $};

\end{tikzpicture}}
\end{center}

By the symmetry property of $\varrho$, we have 
$$\int_0^1 \rho_n (s)\, ds = b-a,$$ 
hence $\rho_n$ is the derivative of a diffeomorphism $g_n$ of $[a,b]$. Since $\rho_n$ is $C^{\infty}$, the diffeomorphism 
$g_n$ is of class $C^{\infty}$. We claim that $g_n$ converges to the identity in the $C^k$ topology for every integer $k$. 
(Hence, by definition, the convergence holds in $C^{\infty}$ topology.) Indeed, for $k \geq 2$, as $n$ goes to infinite, 
we have
$$\| D^k (g_n) \|_{C^0} = \| D^{k-1} (\rho_n) \|_{C^0} = \frac{1}{n} \| D^{k-1} \varrho \|_{C^0} \longrightarrow 0.$$
We will next show that the corresponding diffeomorphism $g$ obtained via (\ref{def-g_n}) is not $C^{1+\alpha}$ for any $\alpha > 0$.

Since $g_n (a) = a$ and $D g_n (a) =1$, we have 
$$Dg (f^n(a)) = D (f^ng_nf^{-n}) (f^n(a)) = \frac{Df^n (g_n (a))}{Df^n (a)} \cdot Dg_n (a) = 1.$$
To compute $Dg (f^n (b))$, first notice that 
$$Df^n (x) = \frac{2^n}{\big[ (2^n-1)x+1 \big]^2}.$$
We compute:  
$$Dg (f^n(b)) 
= \frac{Df^n (g_n(b))}{Df^n (b)}\cdot Dg_n (b) 
= \left[ \frac{(2^n-1)b+1}{(2^n-1)g_n(b)+1} \right] \cdot \left( 1+\frac{1}{n}\right).$$
Since 
$$g_n (b) 
= g_n (a) + \int_a^b \rho_n (s) \, ds 
= g_n (a) + \int_a^b \left( 1 + \frac{\varrho(s)}{n} \right) (s) \, ds 
= a + (b-a) + \int_a^b \frac{\varrho(s)}{n} \, ds
= b + \frac{C}{n},$$
where 
$$I:= \int_a^b \varrho (s) \, ds > 0,$$
this yields
$$Dg (f^n(b)) 
= \left[ \frac{(2^n-1)b + 1}{(2^n-1)(b + I/n) + 1} \right] \cdot \left( 1+\frac{1}{n}\right),$$
hence
$$\left| Dg (f^n (b)) - Dg (f^n (a))\right| 
= \frac{1}{n}  \left[ \frac{(2^n-1)b + 1}{(2^n-1)(b + I/n) + 1} \right] - \left[ \frac{I}{n \, (2^n-1)(b + I/n) + 1} \right].$$
Therefore, for a certain constant $C'$,
\begin{equation}\label{eq-below}
\big| Dg (f^n (b)) - Dg (f^n (a)) \big|  \geq \frac{C'}{n}.
\end{equation}

Finally, putting (\ref{eq-distancia}) and (\ref{eq-below}) together, we obtain
$$\frac{|Dg (f^n(b) - Dg (f^n (a)))|}{|f^n(b)-f^n(a)|^{\alpha}} \geq 
\frac{C'}{C^{\alpha}} \cdot \frac{2^{n\alpha}}{n},$$
which diverges to infinite as $n \to \infty$ provided $\alpha > 0$. This shows that $g$ is not of class $C^{1+\alpha}$.

\vspace{0.1cm}

\begin{rem} It is not hard to modify the preceding example so that the map $f$ has parabolic fixed points. In this 
framework, the divergence in the last step above is much slower, but still holds.
\end{rem}

\vspace{0.35cm}

\noindent{\bf Acknowledgments.} Both authors would like to gratefully acknowledge the guidance of Andr\'es Navas during this work. 
They also thank Crist\'obal Rivas for organizing the seminar on group actions where the problem treated in the paper was presented, 
as well as all the participants of the seminar. 
The first author was funded by the Fondecyt Project 1200114, and the second author by the Fondecyt Project 1181548.


\begin{small}

\vspace{0.3cm}

\noindent Leonardo Dinamarca \& Maximiliano Escayola

\noindent Dpto. de Matem\'atica y Ciencias de la Computación, Univ. de Santiago de Chile (USACH)

\noindent Alameda Lib. Bdo O'Higgins 3363, Estaci\'on Central, Santiago, Chile

\noindent Emails: leonardo.dinamarca@usach.cl, maximiliano.escayola@usach.cl

\end{small}

\end{document}